\documentclass[12pt]{article}
\usepackage{theorem}
\usepackage{amsmath,amssymb,amscd,latexsym}

\theoremstyle{plain}
\theoremheaderfont{\normalfont\bfseries}
\theorembodyfont{\itshape}

\newtheorem{theorem}{Theorem}[section]
\newtheorem{proposition}[theorem]{Proposition}
\newtheorem{lemma}[theorem]{Lemma}

\theorembodyfont{\upshape}
{\theoremstyle{break}

}

\theorembodyfont{\upshape}
\newtheorem{remark}[theorem]{Remark}

\theorembodyfont{\upshape}

\theorembodyfont{\upshape}
\newtheorem{definition}[theorem]{Definition}

\theorembodyfont{\upshape}

\newtheorem{notation}[theorem]{Notation}

\theorembodyfont{\itshape}
\newtheorem{conj}{Conjecture}

\theorembodyfont{\upshape}

\theorembodyfont{\upshape}

\theorembodyfont{\upshape}

\theorembodyfont{\upshape}

\theorembodyfont{\upshape}

\theorembodyfont{\upshape}

\theorembodyfont{\upshape}

\newcommand{\qed}{\nopagebreak\par\hspace*{\fill}$\square$\par\vskip2mm}

\newcommand{\R}{{\mathbb R}}
\newcommand{\Q}{{\mathbb Q}}
\newcommand{\E}{{\cal E}}

\def\sdim{\operatorname{sdim}}
\def\Um{\operatorname{Um}}
\def\MS{\operatorname{MS}}
\def\WMS{\operatorname{WMS}}
\def\GL{\mathop{\mathit GL}}
\def\SL{\mathop{\mathit SL}}
\newcommand{\Aone}{\mathbb A^1}
\def\Row{\operatorname{Row}}
\def\oldpmatrix#1{\begin{pmatrix}#1\end{pmatrix}}

\title{From Mennicke symbols to Euler class groups}
\author{Wilberd  van der Kallen}
\date{}

\begin{document}
\maketitle

\begin{abstract}
Bhatwadekar and Raja Sridharan have constructed a homomorphism
from an orbit set of
unimodular rows to an Euler class group. We show under weaker assumptions
that a generalization of its kernel is a subgroup.
Our tool is a partially defined
operation on the set of unimodular matrices with two rows.
\end{abstract}
\section{The exact sequence}
Consider a commutative noetherian $\Q$-algebra
$A$ of dimension $d$. Let $n=d+1$.
We often assume $n$ is odd.
We try to understand the following exact sequence.
$$0\to \MS_{n-1}(A)\to \Um_{2,n}(A)/E_n(A)\to
\Um_{1,n}(A)/E_n(A)\to \E(A).$$
A good case to keep in mind is $d=6$.
\section{The terms in the sequence}
Let us recall the terms in the sequence.
All matrices will have entries in $A$.
An $m$ by $n$ matrix $M$ with $m\leq n$ is called \emph{unimodular} if
it has a right inverse, which is thus an $n$ by $m$ matrix.
In other words, $M$ is called unimodular if the corresponding map
$A^n\to A^m$ is surjective.
Let $\Um_n(A)=\Um_{1,n}(A)$ be the set of unimodular rows with $n$ entries.
Following Suslin \cite{suslin mennicke}, we say that a Mennicke symbol of order
$n$ on $A$ is
a map $\phi$ from $\Um_n(A)$  to
an abelian group $G$ such that MS1 and MS2 hold:
\begin{itemize}
\item[MS1] For every elementary matrix $\epsilon\in E_n(A)$ and every
$\mathbf{v}\in\Um_n(A)$ we have $\phi(\mathbf{v}\epsilon)=\phi(\mathbf{v})$
\item[MS2] $\phi(x,a_2,\ldots,a_n)\phi(y,a_2,\ldots,a_n)=
\phi(xy,a_2,\ldots,a_n)$
\end{itemize}

Here we
simplified notation from $\phi((x,a_2,\ldots,a_n))$ to
$\phi(x,a_2,\ldots,a_n)$.

\subsection{Universal Mennicke symbols}
The group $\MS_{n-1}(A)$ is by definition the universal target of order
$n-1$ Mennicke symbols on $A$.
We have shown in \cite{vdk module} that $\mathbf{v}$,
$\mathbf{w}\in\Um_{n-1}(A)$
have the same image in $\MS_{n-1}(A)$ if and only if there is
$\epsilon\in E_{n}(A)\cap \GL_{n-1}(A)$ with $\mathbf{v}\epsilon=\mathbf{w}$.
(This needs that $n-1$ is even and $d\leq 2n-7$.)
Moreover, the map $ms:\Um_{n-1}(A)\to \MS_{n-1}(A)$ is surjective, so we
may think of $\MS_{n-1}(A)$ as the orbit set
$\Um_{n-1}(A)/E_{n}(A)\cap \GL_{n-1}(A)$ that has been provided with
a group structure. (For $n-1=2$, $d=1$ this is exactly what Bass Milnor Serre
did for ordinary Mennicke symbols.)

\subsection{The term that we expect to be a group}
The term $\Um_{2,n}(A)/E_n(A)$ is just an orbit set of unimodular
two by $n$ matrices. We expect it also carries a group structure.
Indeed the analogy with algebraic topology, cf.\ section \ref{Banach},
predicts an abelian group
structure for $d\leq2n-6$. (Such an analogy gave correct predictions
in \cite{vdk module}.)
But for now we have no such group structure, so \emph{our
exact sequence will be
one of pointed sets only!}

The map $\MS_{n-1}(A)\to \Um_{2,n}(A)/E_n(A)$ we define by sending
$ms(\mathbf{v})$ to the orbit of $\oldpmatrix{1&0\cr0&\mathbf{v}}$.
So it is some kind of stabilization map, and we are grasping at an
unstable $K$-theory in which non-square unimodular matrices are employed.
In other words, we are exploring the connection between $K$-theory and
Grassmannians, albeit in an idiosyncratic manner.

{}From what we said above one sees that $\MS_{n-1}(A)\to \Um_{2,n}(A)/E_n(A)$
 is well defined. One also sees that the sequence
is exact at $\MS_{n-1}(A)$. The map $\Row_1:\Um_{2,n}(A)/E_n(A)\to
\Um_{1,n}(A)/E_n(A)$ associates to the orbit of a matrix the orbit of its
first row. Exactness of our sequence at $\Um_{2,n}(A)/E_n(A)$ is then easy.
\subsection{The orbit group}
Recall that we have shown in \cite{vdk group} that
$\Um_{1,n}(A)/E_n(A)$ carries a group structure for $n=d+1\geq3$.
This was extended in \cite{vdk module} to the range $2\leq d\leq 2n-4$.

\subsection{The Euler class group $\E(A)$}
The last term in our sequence is, as our title indicates, the Euler class
group $\E(A)$ of Bhatwadekar and Raja Sridharan as introduced in
\cite{BhatwadekarSridharan} for
a commutative noetherian $\Q$-algebra $A$ of dimension $d$.
Let us recall what it is about.
We fix a rank one projective $A$-module $L$. We think of $L$ as an
orientation bundle. An $L$-oriented rank $d$ bundle is a projective $A$-module
$P$ of rank $d$ together with an orientation isomorphism $\det(P)\to L$.
A generic section of such a bundle is an $A$-module map $P\to A$ whose
image $J$ is an ideal of height at least $d$. So $A/J$ has \emph{finite}
spectrum,
and $J=A$ is allowed. We say the section vanishes at the maximal ideal
$\mathfrak m$ if  $\mathfrak m$ contains the \emph{ideal of vanishing}
$J$. By taking a formal sum of the zeroes of a section one would get
a zero cycle,
which one could then consider modulo rational equivalence.
That is not quite what we will need, however.
The zeroes need to be enriched with
orientation data (and some form of multiplicity).
Thus we now want to say what an
\emph{oriented local zero} of a section is.
Basically one gets it by localizing and completing at
$\mathfrak m$,  when $\mathfrak m$ contains
$J$. More precisely, one tensors the orientation $\det(P)\to L$
and the section $P\to J$ with $A_{\mathfrak m}/JA_{\mathfrak m}$.
That gives a representative of
an oriented local zero. Another representative, coming from a surjective map
$Q\to JA_{\mathfrak m}/J^2A_{\mathfrak m}$,
where $Q$ is an $L\otimes A_{\mathfrak m}/JA_{\mathfrak m}$-oriented
projective
$A_{\mathfrak m}/JA_{\mathfrak m}$-module of rank $d$,
defines the same oriented local zero at $\mathfrak m$ if there is an
isomorphism $Q\to P\otimes A_{\mathfrak m}/JA_{\mathfrak m}$
that is compatible with both
the sections and the orientations. So for an oriented local zero we need
an $\mathfrak m$, a height $d$ ideal $J_{\mathfrak m}$ in $A_{\mathfrak m}$,
an $L\otimes A_{\mathfrak m}/J_{\mathfrak m}$-oriented projective
$A_{\mathfrak m}/J_{\mathfrak m}$-module of rank $d$,
and finally a surjection from that
projective module onto $JA_{\mathfrak m}/J^2A_{\mathfrak m}$.
One now defines the Euler class group $\E(A,L)$ to be
an abelian group whose generators are oriented local zeroes
(sums of them are also known
as oriented zero cycles) and
whose relations say that
 any two
generic sections of the same $L$-oriented bundle
yield the same sum total of oriented local zeroes.
This common total is called the Euler class of the $L$-oriented bundle and
the big theorem in \cite{BhatwadekarSridharan} says that $P$ has a
nowhere vanishing global section exactly when its Euler class vanishes.
Actually they define $\E(A,L)$ with fewer relations, and then they show that
the relations
we just described follow from their set.

Now take for $L$ the trivial line bundle $A$ and put $\E(A)=\E(A,A)$.
How do we define the map $\Um_{1,n}(A)/E_n(A)\to \E(A)$?
That is easy. If $\mathbf{v}\in\Um_{1,n}(A)$, then the kernel $P$
of $\mathbf{v}: A^n\to A$ is naturally an $A$-oriented projective of
rank $d$.
We just send the orbit of  $\mathbf{v}$ to the Euler class of $P$.

\section{The homomorphism $\Um_{1,n}(A)/E_n(A)\to\E(A)$}
It was discovered by Bhatwadekar and Raja Sridharan that this
defines a homomorphism!
They showed it first for the case $d=2=n-1$ pioneered by Vaserstein.
In that case one actually understands how
$\Um_{1,n}(A)/E_n(A)$ and $\E(A)$
are related. Vaserstein constructed the group structure on $\Um_{3}(A)/E_3(A)$
by using that an $A$-oriented projective module
of rank two is just a projective module
of rank two with a nondegenerate alternating bilinear form on it.
Then he used the group structure on the appropriate Witt group of alternating
forms.
His group structure thus had a meaning. In \cite{vdk group} we showed that
the group structure on $\Um_{n}(A)/E_n(A)$ for $d=n-1\geq2$
could be constructed by induction on $d$. Following this induction
Bhatwadekar and Raja Sridharan then establish that one gets a homomorphism
$\Um_{d+1}(A)/E_{d+1}(A)\to \E(A)$ for larger $d$ also.
Actually this induction is rather funny: Recall that Bass has
shown that for even $n$ any unimodular row of
length $n$ over a commutative ring
 is the first row of a unimodular two by $n$ matrix
\cite[I.4.12]{lam}, or that the
corresponding projective has a free summand. Therefore for even $n$
the homomorphism $\Um_{1,n}(A)/E_n(A)\to\E(A)$ is trivial and the induction
has to increase $d$ in steps of size two.

An intriguing byproduct of the existence of this
homomorphism is that one gets a \emph{subgroup} by taking
 those unimodular rows of
length $n=d+1$
which define
a projective module with nowhere vanishing global section (in other words,
with a free summand).
This will lead to curious facts like:
If the kernel of $(a_1,a_2,\ldots,a_n)$ has a free summand, so does
the kernel of $(a_1^7,a_2,\ldots,a_n)$. Does this hold only under conditions
on $d$? (Note that we are using  the condition $n\geq d+1$ here,
while in \ref{subgroup}
it will be $d\leq 2n-5$.)

We want to understand the homomorphism and the subgroup in other ways.
First of all, we have in  \cite{vdk module}
a more mysterious construction of the group structure, by generators and
relations.
This more mysterious construction works as soon as  $2\leq d\leq 2n-4$.
So one would expect to construct a homomorphism
$\Um_{d+1}(A)/E_{d+1}(A)\to \E(A)$
by checking that the relations hold in $\E(A)$.

When trying to do this, one finds that the defining relations used
in \cite{vdk module} are not so obvious in $\E(A)$.
In a similar situation I had some help from Ofer Gabber who showed me
that in Borsuk's cohomotopy groups one has relation MS3:

\begin{itemize}
\item[MS3] $\phi(x,a_2,\ldots,a_n)\phi(y,a_2,\ldots,a_n)=
\phi(xy,a_2,\ldots,a_n)$ if $x+y=1$.
\end{itemize}

Now that we are at it,
let us list a few more relations that may or may not hold.

\begin{itemize}
\item[MS4] $\phi(f^2,a_2,\ldots,a_n)\phi(g,a_2,\ldots,a_n)=
\phi(f^2g,a_2,\ldots,a_n)$
\item[MS5]$\phi(r,a_2,\ldots,a_n)
\phi(1+q,a_2,\ldots,a_n) = \phi(q,a_2,\ldots,a_n)$,
if $r(1+q)\equiv q\bmod (a_2,\ldots,a_n)$
\item[MS6] $\phi(x,a_2,\ldots,a_n)=
\phi(-x,a_2,\ldots,a_n)$
\item[MS7] $\phi(x,a_2,\ldots,a_n)^m=
\phi(x^m,a_2,\ldots,a_n)$ for $m\geq2$
\end{itemize}

We showed in \cite{vdk module} that for $2\leq d\leq 2n-4$ the universal
MS1\&MS5 symbol defines a bijection $wms$ from $\Um_n(A)/E_n(A)$
to an abelian group
$\WMS_n(A)$. We also saw there that $wms$ satisfies MS4.
Despite Gabber's observation we did not suspect that $wms$ satisfies
MS3. This had to wait till the work of
Bhatwadekar and Raja Sridharan. In fact one has the simple

\begin{lemma}
Under \textup{MS1\&MS4} relations \textup{MS3} and \textup{MS5} are equivalent.
\end{lemma}
\paragraph{Proof}
Try $\overline{(1+q)y}\equiv \overline1\bmod (a_2,\ldots,a_n)$.\qed

\

Now recall the following variation of the Mennicke-Newman lemma

\begin{lemma}\label{tails}
Let $d\leq 2n-3$. Let $\mathbf{v}$, $\mathbf{w}\in\Um_{n}(A)$.
There are $\epsilon$, $\delta\in E_n(A)$ and $x$, $y$, $a_i\in A$
so that $\mathbf{v}\epsilon=(x,a_2,\ldots,a_n)$,
$\mathbf{w}\delta=(y,a_2,\ldots,a_n)$, $x+y=1$.
\end{lemma}
\paragraph{Proof} (Backwards)
Once we have $v_1+w_1$ equal to one, it is easy to get $v_i-w_i$ zero for $i>1$
by adding multiples of $v_1$ to $v_i$ and of $w_1$ to $w_i$.
We may change $v_1+w_1$ by adding multiples of
$v_i$ to $v_1$ and of $w_i$ to $w_1$.
So we would be through if $(v_2,\ldots,v_n,w_1,\ldots,w_n)$
were unimodular.
Now observe that $(v_1w_1,v_2,\ldots,v_n,w_1,\ldots,w_n)$
is unimodular because $A$ is commutative. (This may be checked over residue
fields.)
By a stable range condition we may add multiples of $v_1w_1$
to the other entries of $(v_1w_1,v_2,\ldots,v_n,w_1,\ldots,w_n)$
to achieve $(v_2,\ldots,v_n,w_1,\ldots,w_n)$
is indeed unimodular.
(It was Keith Dennis who told me about this trick of Vaserstein
to treat two rows together in such manner.
The trick is crucial for obtaining all the results
here and in \cite{vdk module} that have a factor 2 in
front of the $n$.)
\qed

\

We can now change the description by generators and relations of the
$\WMS_n(A)$ groups.

\begin{theorem}
For $2\leq d\leq 2n-4$ the universal
\textup{MS1\&MS3} symbol defines a bijection $wms$ from $\Um_n(A)/E_n(A)$ to an
abelian
group
$\WMS_n(A)$.
\end{theorem}

\begin{remark}
Note that it is rather obvious from lemma~\ref{tails}
that the target of the
universal MS1\&MS3 symbol has to be abelian. In \cite{vdk module} we
had much more trouble to understand why $\WMS_n(A)$ is abelian.
\end{remark}

\paragraph{Proof of Theorem}
We compare the universal MS1\&MS3 symbol $\phi:  \Um_n(A)\to G$ with the
universal
MS1\&MS5 symbol $wms:\Um_n(A)\to\WMS_n(A)$
of \cite{vdk module}.
Clearly we have a homomorphism $\tau: G\to\WMS_n(A)$ so that
$wms=\tau\circ\phi$. An element of $\ker\tau$ may be written as
$\phi(\mathbf{v})-\phi(\mathbf{w})$. As $wms$ defines a bijection
from $\Um_n(A)/E_n(A)$ to $\WMS_n(A)$, the orbits of $\mathbf{v}$ and
$\mathbf{w}$ must be the same. But then by MS1 the element
$\phi(\mathbf{v})-\phi(\mathbf{w})$ is zero. As $wms$
is also surjective, $\tau$ is an isomorphism.\qed

\

Now we are ready to give another construction of the
homomorphism $\Um_{1,n}(A)/E_n(A)\to\E(A)$ for $d+1=n\geq 3$.
All we have to check is that the relation MS3 is satisfied.
We do not care whether $n$ is odd or even.

Suppose we are given unimodular
$(x,a_2,\ldots,a_n)$ and $(y,a_2,\ldots,a_n)$ with
 $x+y=1$. By modifying these representatives of $E_n(A)$-orbits we
may assume $(x,a_2,\ldots,a_{n-1})$ and $(y,a_2,\ldots,a_{n-1})$ both
generate
an ideal of height at least $d$. (Recall $A$ contains $\Q$ if you
do not wish to consider residue fields with two elements in this prime
avoidance exercise.)

Say $\mathbf{v}=(v_1,\ldots,v_n)$ is unimodular so that the ideal
$J=(v_1,\ldots,v_{n-1})$ has height at least $d=n-1$ and $P$
is the corresponding projective module. One may map $P$
to $A$ with image $J$ by sending every vector in $A^n$ to its last coordinate
and then restricting that map to $P$.
(That the image is $J$ follows from the exactness of the Koszul
complex of $\mathbf{v}$.)
So that defines a suitable generic section.
The oriented local zero at
${\mathfrak m}$ of this section may be represented by the map
$(A_{\mathfrak m}/JA_{\mathfrak m})^{n-1}\to
JA_{\mathfrak m}/J^2A_{\mathfrak m}$
given by the matrix $(v_1,\ldots,v_{n-1})$, where the free
module $(A_{\mathfrak m}/JA_{\mathfrak m})^{n-1}$
has its standard orientation multiplied by a possible sign and a
certain power of
the unit
$v_n\bmod JA_{\mathfrak m}$. (Exercise, compare \cite{BhatwadekarSridharan}.)

Apply this for $\mathbf{v}=(x,a_2,\ldots,a_n)$, $(y,a_2,\ldots,a_n)$,
$(xy,a_2,\ldots,a_n)$ respectively.
At a maximal ideal ${\mathfrak m}$ that contains the ideal
$J_x=(x,a_2,\ldots,a_{n-1})$
we have $y\equiv 1\bmod J_x$, so the map
$(A_{\mathfrak m}/J_xA_{\mathfrak m})^{n-1}\to
J_xA_{\mathfrak m}/J_x^2A_{\mathfrak m}$
is identical to its analogue
$(A_{\mathfrak m}/J_{xy}A_{\mathfrak m})^{n-1}\to
J_{xy}A_{\mathfrak m}/J_{xy}^2A_{\mathfrak m}$.
This way we see that the total of the local contributions for
$(xy,a_2,\ldots,a_n)$ is indeed the sum of the two totals for
$(x,a_2,\ldots,a_n)$ and $(y,a_2,\ldots,a_n)$ respectively.
\qed

\section{The image of $\Row_1$}
Bhatwadekar and Raja Sridharan pointed out that
the kernel of the
map  $\Um_n(A)/E_n(A)\to \E(A)$ gives an interesting
subgroup of $\WMS_n(A)=\Um_n(A)/E_n(A)$.
(Interesting subgroups could still be trivial.)
We will now approach this subgroup from the other side, viewing it
as the image of $\Row_1:\Um_{2,n}(A)/E_n(A)\to
\Um_{1,n}(A)/E_n(A)$. While the Euler class homomorphism only works for
$n=d+1$, we will see that in fact the image of $\Row_1$ is a subgroup
for $d\leq2n-5$. (Bhatwadekar and Raja Sridharan inform me that if
$A$ is a regular ring containing a field,
one could
push their construction to obtain their homomorphism in a similar wider range,
adapting the target group to this purpose.)
We will show that the image of $\Row_1$ is a subgroup
by showing that in some sense
$\Row_1$ is a homomorphism. Recall that we do not have a
group structure yet on $\Um_{2,n}(A)/E_n(A)$.
We will define an  `operation' $*$ on
$\Um_{2,n}(A)/E_n(A)$
that is probably well-defined. As our solution of the analogous
problem in \cite{vdk module} is tedious and not illuminating,
we have not attempted
to prove directly that $*$ is well defined.
We conjecture that $*$ describes the operation in a group
structure. In any case,
we will show that  $\Row_1$ is a `homomorphism'
and conclude from this that the image of  $\Row_1$ is a subgroup.
In other words, just like in \cite{vdk module}
our proofs are proofs, but they leave room for improvement.
As $\Row_1$ is surjective if $n$ is even, we are mostly interested in odd
$n$.

First we need a better understanding of $\Um_{2,n}(A)/E_n(A)$.
\begin{lemma}\label{left}
Let $A$ be commutative, $m< n$ and $n\geq3$.
The action by left multiplication
of $E_m(A)$ on $\Um_{m,n}(A)/E_n(A)$
is trivial.
\end{lemma}
\paragraph{Proof}
There are simple
explicit formulas as in \ref{over} below,
but the proof I like most uses the
following reasoning. We wish to show that an elementary
matrix $e_{ij}(t)$ acts trivially
on the class of $M\in \Um_{m,n}(A)$.
But then one may as well view $t$ as a variable and recall there are
Quillen local global principles. We use the formalization by Vaserstein as
discussed in \cite{lam}, combined with Suslin's proof \cite{suslin pol}
that a relative elementary
group is normal. It shows that we may assume $t$ is the variable in
a polynomial ring $A$ over a local ring which contains the entries of $M$.
Then we may first transform
$M$ to standard form using column operations.
\qed

\begin{lemma}\label{adjoint orbit}
Let $A$ be commutative, $m< n$ and $n\geq3$.
If $M_1$, $M_2\in \Um_{m,n}(A)$ have the same right inverse $N$, then they
belong to the same orbit under $E_n(A)$.
\end{lemma}
\paragraph{Proof}
For any $t$ the matrix $N$ is a right inverse of $tM_1+(1-t)M_2$.
We  may again pretend $t$ is a variable, apply a local global principle,
and reduce to the case that $N$ consists of the first $m$ columns of the
identity matrix.
\qed

\begin{remark}
Obviously these two lemmas are still valid for $m=n$.
\end{remark}

\begin{remark}
Applying this kind of reasoning we also see that for $m<n$, \hbox{$n\geq3$},
our orbit sets $\Um_{m,n}(A)/E_n(A)$
of unimodular matrices may be identified
with orbit sets of \emph{split unimodular matrices}, meaning pairs
consisting of a unimodular
matrix and a right inverse. That may be the better way of looking
at it.
\end{remark}

We need a peculiar analogue of the
Mennicke Newman lemma in $\Um_{2,n}(A)$ for $n\geq3$.

\begin{lemma}\label{prod}
Let $n\geq4$ and $d\leq2n-5$.
Let two elements of $\Um_{2,n}(A)/E_n(A)$ be given.

Then we may choose representatives of the form
$\oldpmatrix{a&b&y_{11}&y_{12}&z_{1,1}\cdots z_{1,n-4}\cr
g&-a&y_{21}&y_{22}&z_{2,1}\cdots z_{2,n-4}}$ and
$\oldpmatrix{1-a&-b&y_{11}&y_{12}&z_{1,1}\cdots z_{1,n-4}\cr
-g&1+a&y_{21}&y_{22}&z_{2,1}\cdots z_{2,n-4}}$ respectively.
With such choice put $$X=\oldpmatrix{a&b\cr g&-a},\quad
Y=\oldpmatrix{y_{11}&y_{12}\cr y_{21}&y_{22}},\quad
Z=\oldpmatrix{z_{1,1}\cdots z_{1,n-4}\cr
z_{2,1}\cdots z_{2,n-4}},\quad I=\oldpmatrix{1&0\cr 0&1}
.$$ Note that $X$ has trace zero.

Then the matrix
$\oldpmatrix{X(I-X)&(I-X)Y+YX&Z}$
is also unimodular.
\end{lemma}
\paragraph{Proof}
For simplicity of notation we put $n=5$.
We know from~\ref{tails} how to get the $(1,1)$ entries to sum to one.
Next use row operations to get the first column into shape.
So now we have representatives
$\oldpmatrix{a&b&c&e&f\cr g&z&h&j&k}$ and $\oldpmatrix{1-a&r&s&t&u\cr
-g&v&w&x&y}$.
Use column operations to add multiples of $g$
to $h$, $j$, $k$, $w$, $x$, $y$ so as to get
 $\oldpmatrix{vz&h&j&k&w&x&y}$ unimodular (\cite{vas rank}).
Then use more column operations to get that
$\oldpmatrix{h&j&k&w&x&y}$ is unimodular.
Then use column operations to get $z+v$ equal to one. And add a
$(-r-b)$-multiple of
the first column to the second in both representatives.
That makes  $b+r$
zero. Put $X=\oldpmatrix{a&b\cr g&z}$. Now the first representative starts
with $X$ and the second with $I-X$. So with column operations we can get
the later columns equal. Of course this may destroy
unimodularity of $\oldpmatrix{h&j&k&w&x&y}$.
Our representatives look like
$\oldpmatrix{a&b&c&e&f\cr g&z&h&j&k}$ and
$\oldpmatrix{1-a&-b&c&e&f\cr -g&1-z&h&j&k}$ by now.
{}From now one we will keep them in this form.
By column operations we can add any multiple of $\det(X)\det(I-X)$ to
any of
$c$, $e$, $f$,  $h$, $j$, $k$.
So we can make
$\oldpmatrix{c&e&f&h&j&k}$ unimodular.
Adding multiples of the later columns to $X$, we can then further achieve
that $X$ has trace zero. The representatives are now in the required form.
Over a residue field of $A$ the column span of
$X-X^2=\oldpmatrix{a-a^2-bg&b\cr g&-a-a^2-bg}$ equals the column span of
$X$ or of $I-X$, or both.
Therefore  $\oldpmatrix{X-X^2&Y&Z}$
is also unimodular.
But we want that
$\oldpmatrix{X(I-X)&(I-X)Y+YX&Z}$ is unimodular.
We only need to look at residue fields where $X(I-X)$ is singular. So we
may assume that $a^2+bg$ equals zero or one.
If $a^2+bg$ is zero, then $X^2=0$ and the column span of
$\oldpmatrix{X(I-X)&(I-X)Y+YX&Z}$ is that of
$\oldpmatrix{X(I-X)&Y(I+X)&Z}$, or of
$\oldpmatrix{X&Y(I+X)(I-X)&Z}=\oldpmatrix{X&Y&Z}$.
If $a^2+bg$ is one, then $X^2=I$ and the column span of
$\oldpmatrix{X(I-X)&(I-X)Y+YX&Z}$ is the same as the
column span of
$\oldpmatrix{(I-X)X&YX&Z}$, or of
$\oldpmatrix{I-X&YX^2&Z}=\oldpmatrix{I-X&Y&Z}$.\qed

\begin{notation}\label{star}
If $[M]$, $[N]$, $[T]$ in $\Um_{2,n}(A)/E_n(A)$
 have respective representatives of the form
$\oldpmatrix{X&Y&Z}$,
 $\oldpmatrix{I-X&Y&Z}$,
$\oldpmatrix{X(I-X)&(I-X)Y+YX&Z}$ with $I$,  $X$,
$Y$, $Z$ as in the lemma, then
we write $$[T]=[M]*[N].$$
The notation is just a notation. It does not mean that we
claim $*$ is a well defined
operation on $\Um_{2,n}(A)/E_n(A)$,
even though we expect this to be the case.
All we really have is a partially defined
operation on $\Um_{2,n}(A)$.
\end{notation}

\begin{conj}
The operation $*$ is well defined and
gives the expected group operation on $\Um_{2,n}(A)/E_n(A)$ for
$n\geq4$ and
$d\leq2n-6$.
\end{conj}

\begin{remark}
 Our formula is the correct
one after inverting $\det(X)$ or $\det(I-X)$.
Indeed
the $E_n(A[\det(X)^{-1}])$-orbit of $\oldpmatrix{X&Y&Z}$
is the trivial orbit and the $E_n(A[\det(X)^{-1}])$-orbit of
$\oldpmatrix{X(I-X)&(I-X)Y+YX&Z}$
contains $\oldpmatrix{I-X&Y&Z}$.
Similarly the $E_n(A[(1+\det(X))^{-1}])$-orbit of
$\oldpmatrix{I-X&Y&Z}$ is the trivial orbit
and the $E_n(A[(1+\det(X))^{-1}])$-orbit of
$\oldpmatrix{X(I-X)&(I-X)Y+YX&Z}$
contains $\oldpmatrix{X&Y&Z}$.\end{remark}

\begin{lemma}
Let $n\geq4$ and $d\leq2n-5$.
If $[T]=[M]*[N]$ then $$\Row_1([T])=\Row_1([M])\Row_1([N]).$$
\end{lemma}
\paragraph{Proof}Using that MS1 holds in
the group $\Um_n(A)/E_n(A)$,
one sees that
$\Row_1([M])$ equals $\Row_1(\oldpmatrix{X&(I-X)Y+YX&Z})$.
Similarly one sees
that $\Row_1([N])$ is the same as $\Row_1\oldpmatrix{(I-X)
\oldpmatrix{1&0\cr0&-1}&
(I-X)Y+YX&Z}$. Now use that MS3 holds.
\qed

\


We conclude that the image of $\Row_1$ is closed under taking products.
To deal with inverses we first extend the
 Swan--Towber Krusemeyer Suslin lemma that a unimodular row
of the form $(a^2,b,c)$ is completable to an invertible three by three
matrix.

\begin{lemma}
Let $A$ be commutative, $n\geq3$, $\mathbf{a}=
(a_1^2,a_2,\ldots, a_n)\in \Um_n(A)$.
Then $\mathbf{a}$ is in the image of $\Row_1:\Um_{2,n}(A)\to\Um_n(A)$.
\end{lemma}
\paragraph{Proof}
Indeed by Bass one does not need the square if $n$ is even, and if $n$ is
odd one argues as follows.
For $n=3$ it is the Swan--Towber Krusemeyer Suslin lemma \cite{Krusemeyer}.
Now let $J$ be the ideal
generated by $a_4,\ldots,a_n$. Complete $(a_1^2,a_2,a_3)$ modulo $J$
to a three by three
matrix $\bar{M}$ of determinant one, lift to a matrix $M=(m_{ij})$
with first row $(a_1^2,a_2, a_3)$,
and choose $b_4,\ldots,b_n$ so that
$\Sigma_{i=2}^{(n-1)/2}(a_{2i}b_{2i+1}-a_{2i+1}b_{2i})=\det(M)-1$.
Then the following matrix is in $\Um_{2,n}(A)$:
$\oldpmatrix{a_1^2&a_2& a_3&a_4&\cdots&a_n\cr
m_{21}&m_{22}&m_{23}&b_4&\cdots&b_n}$.
\qed

\begin{remark}
Bhatwadekar and Raja Sridharan also show \cite[5.4]{BhatwadekarSridharan}
that for Euler classes multiplication
by a square does not make a difference.
\end{remark}

We are now ready to prove

\begin{theorem}\label{subgroup}
Let $A$ be a commutative noetherian ring of Krull
dimension $d$. (Or, more generally, a commutative ring whose
stable range dimension $\sdim$ \cite{vdk module} equals $d$.)
If $d\leq2n-5$ then the image of
$\Row_1:\Um_{2,n}(A)/E_n(A)\to
\Um_{n}(A)/E_n(A)$ is a subgroup.
\end{theorem}
\paragraph{Proof}
We know already that the image is closed under product.
If $\mathbf{a}=
(a_1,a_2,\ldots, a_n)\in \Um_n(A)$  has as right inverse the  transpose of
$(b_1,b_2,\ldots, b_n)$, then the inverse of $wms(\mathbf{a})$
equals the product $wms(-a_1,a_2,\ldots, a_n)wms(b_1^2,a_2,\ldots, a_n)$
by \cite[3.5]{vdk module}.\qed

\

For $n$ even $\Row_1$ is surjective, but for odd $n$ it hits only very special
elements. We now show that
these elements satisfy Rao's relation MS6. Recall that Rao used MS6
in \cite{rao three} to show that elements are trivial.

\begin{lemma}Let $n$ be odd, $d\leq2n-4$.
Rao's relation \textup{MS6} holds in the image of $\Row_1$
and therefore \textup{MS7} also holds in this image.
\end{lemma}
\paragraph{Proof}That MS6 implies MS7 is proved in \cite{rao three}.
To prove MS6 we consider $M\in \Um_{2,n}(A)$ with right inverse
$N$. Say $\mathbf{v}$ is the first row of $M$, $\mathbf{w}$
is the second row, and $\mathbf{z}$
is the second column of $N$. Let $D$ be the $n$ by $n$
diagonal matrix with diagonal $(-1,1,\ldots,1)$.
By Vaserstein's Whitehead lemma the matrix
$S=(I_n-2\mathbf{z}\mathbf{w})D$ is in $E_{n+1}(A)\cap \GL_n(A)$,
so $mse(\mathbf{v})=mse(\mathbf{v}S)$ by \cite[5.3]{vdk module}.
But $\mathbf{v}S=\mathbf{v}D$.\qed

\begin{remark}
In fact, if $t$ is a unit, one similarly has $\phi(x,a_2,\ldots,a_n)=
\phi(tx,a_2,\ldots,a_n)$ for $(x,a_2,\ldots,a_n)$ in the image of $\Row_1$.
\end{remark}

\section{$\Aone$-homotopy}\label{A1hom}
\begin{definition}
We say that $\bf v$, ${\bf w}\in \Um_n(A)$ are $\Aone$-homotopic,
notation ${\bf v}\sim_{\Aone}\bf w$, if there is
${\bf z}(t)\in \Um_n(A[t])$ with ${\bf z}(0)=\bf v$, ${\bf z}(1)=\bf w$.
\end{definition}

\begin{lemma}
If ${\bf v}\sim_{\Aone}\bf w$ and $\epsilon\in E_n(A)$ then
${\bf v}\sim_{\Aone}\bf w\epsilon$.
\end{lemma}
\paragraph{Proof}
 Choose $\alpha(t)\in E_n(A[t])$ with $\alpha(0)=1$ and
$\alpha(1)=\epsilon$. Then multiply the homotopy by $\alpha(t)$.\qed

\begin{lemma}
Let $d\leq2n-4$.
The orbits $[{\bf v}]:={\bf v}E_n(A)$
 with  ${\bf v}\sim_{\Aone}(1,0,\ldots,0)$ form a subgroup $H$ of
$\Um_{n}(A)/E_n(A)$. Its cosets are the equivalence classes of the relation
$\sim_{\Aone}$.
\end{lemma}
\paragraph{Proof}
The dimension of $A[t]$ is $d+1$.
To see that $H$ is closed under product, use
lemma~\ref{tails}, MS3 and the previous lemma.
To deal with inverses, compare with the proof of theorem~\ref{subgroup}.
Arguing in the same manner one finds that ${\bf v}\sim_{\Aone}\bf w$
if and only if $[{\bf v}]^{-1}[{\bf w}]\in H$. Thus $\sim_{\Aone}$ is an
equivalence relation and its classes are the cosets of $H$.
\qed

\begin{notation}
We write $\Um_n(A)/\sim_{\Aone}$ for the set of $\Aone$-homotopy classes in
$\Um_n(A)$.
\end{notation}

\begin{proposition}
Let $d\leq2n-4$.
The group structure on $\Um_{n}(A)/E_n(A)$  induces
a group structure on $\Um_n(A)/\sim_{\Aone}$.\qed
\end{proposition}

\begin{remark}
Similarly one could put an equivalence relation on $\Um_{2,n}(A)$
generated by $\Aone$-homotopy. The operation $*$ in \ref{star}
should describe
a group structure on the $\Aone$-homotopy equivalence classes.
\end{remark}

\section{Orbit sets over Banach algebras}\label{Banach}
Let us now take $A$ to be the Banach algebra of continuous
real valued functions on some
finite $d$-dimensional
CW complex $X$. Then one knows \cite{vas mennicke}
 that for $n\geq3$ the orbit set
$\Um_n(A)/E_n(A)$ is in bijective correspondence with the set $[X,\R^n-0]$
of homotopy classes of maps from $X$ to $\R^n-0=\Um_n(\R)$.
This gives a topological explanation why for $2\leq d\leq2n-4$
one has a group structure on $\Um_n(A)/E_n(A)$.
We now want to show that similarly
$\Um_{m,n}(A)/E_n(A)$ is in bijective correspondence with
the set $[X,\Um_{m,n}(\R)]$ if $m<n$, $n\geq3$.
Then for $d\leq 2n-6$
the fact that the image of
$\Row_1:\Um_{2,n}(A)/E_n(A)\to
\Um_{1,n}(A)/E_n(A)$ is a subgroup may be
explained by the fact that the stable category is an additive category.

\begin{notation}
If $M$, $N$ are matrices, let $M\perp N$ denote their block sum
$\oldpmatrix{M&0\cr 0&N}$.
\end{notation}

\begin{lemma}\label{over}
Let $m<n$ and $n\geq3$. Let $T\in GL_m(A)$ and let $M\in \Um_{m,n}(A)$.
Then $TM$ and $M(T\perp I_{n-m})$ are in the same orbit under
$\GL_n(A)\cap E_{m+n}(A)$.
If $T$ is close to the identity then they are in the same $E_n(A)$-orbit.
\end{lemma}
\paragraph{Proof}Let $N$ be a right inverse of $M$.
By Vaserstein's Whitehead lemma $I_n+N(T-I_m)M$ and
$T=I_m+(T-I_m)MN$ are in the same coset of $E_{m+n}(A)$.
This shows that $TM=M(I_n+N(T-I_m)M)$ and $M(T\perp I_{n-m})$
are in the same orbit under
$\GL_n(A)\cap E_{m+n}(A)$.
If $T$ is close to one, both $T$ and $I_n+N(T-I_m)M$ are close to the
identity. Now recall that $\SL_n(A)\cap E_n(A)$ is open.
\qed

\begin{lemma}
Let $m<n$ and $n\geq3$.
The action of $E_n(A)$ on $\Um_{m,n}(A)$ has open orbits.
\end{lemma}
\paragraph{Proof}
Let $\GL_n^0$ be the path connected component of the identity in the group
$\GL_n(A)$. The kernel of the determinant map $\GL_n^0\to A^*$
is contained in $E_n(A)$, so
the action of $\GL_n^0$ on $\Um_{m,n}(A)/E_n(A)$ factors
through the determinant.
But the previous lemma tells us that if $M\in \Um_{m,n}(A)$
and $t$ is close to one in $A$,  then $tM$ is in the same orbit as
$M(tI_m\perp I_{n-m})$.
So $M$ is in the same orbit as $M(I_m\perp t^{-1}I_{n-m})$.
The determinant of $I_m\perp t^{-1}I_{n-m}$ covers a neighborhood of 1 as
$t$ varies. So the stabilizer in $\GL_n^0$ of the
$E_n(A)$-orbit of $M$ is open,
and must be all of $\GL_n^0$, as $\GL_n^0$ is connected. Now if $N$ is
a right inverse of $M$ and $M'$ is close to $M$, then $M'N$ is in
the path connected component of the identity in $\GL_m(A)$.
Now use lemma \ref{over} and lemma \ref{adjoint orbit}.
\qed

\begin{remark} Let $n\geq4$ and $d\leq2n-6$. The operation $*$ in \ref{star}
was selected
so that it agrees with
the group structure on $\Um_{2,n}(A)/E_n(A)$, as described by Borsuk
\cite{Borsuk}.
\end{remark}

\end{document}